\documentclass{amsart}

\usepackage{amsfonts}
\usepackage{amssymb}
\usepackage{amsmath}
\usepackage{hyperref}
\usepackage{graphicx}

\newcommand{\Gam}{\Gamma_{\mathbb{R}}}
\newcommand{\half}{\frac{1}{2}}
\newcommand{\thalf}{\tfrac{1}{2}}
\newcommand{\summ}{\mathop{{\sum}^{\star}}}

\newcommand{\sums}{\mathop{{\sum}^{*}}}

\newcommand{\sym}{\text{sym}}
\newcommand{\intt}{\int_{-\infty}^{\infty}}

\newcommand{\lam}{\lambda_f}

\numberwithin{equation}{section}

\newtheorem{theorem}{Theorem}[section]

\newtheorem{lemma}[theorem]{Lemma}

\makeatletter
\@namedef{subjclassname@2010}{%
  \textup{2010} Mathematics Subject Classification}
\makeatother

\begin{document}

\title{A mean value of a triple product of $L$-functions }

\author{Jack Buttcane}

\thanks{J. Buttcane: Department of Mathematics, University at Buffalo, SUNY, Buffalo, NY 14260; email: buttcane@buffalo.de}

\author{Rizwanur Khan}

\thanks{R. Khan (corresponding author): Science Program, Texas A\&M University at Qatar, PO Box 23874, Doha, Qatar; email: rk2357@gmail.com}



\begin{abstract} 
Luo has proven an optimal upper bound for the $L^4$-norm of dihedral Maass forms of large eigenvalue, by bounding a mean value of triple product $L$-functions. Motivated by this result, we study a mean value of $L$-functions having similar shape, and obtain for it an asymptotic with power savings. Our work may be helpful in eventually obtaining an asymptotic for the $L^4$-norm.
\end{abstract}

\maketitle

\section{Introduction}

This paper is motivated by a problem in arithmetic quantum chaos, which we describe first. The so called random wave conjecture \cite{ber,hejrac, hejstr} states for $ \Gamma_0(d) \backslash \mathbb{H}$ that any Hecke-Maass cusp form $f$ with large Laplacian eigenvalue $\lambda_f $ should have Gaussian moments (and therefore behave like a random wave). More precisely in the case of the fourth moment, it is conjectured that with the normalization
\begin{align}
\frac{1}{\int_{ \Gamma_0(d) \backslash \mathbb{H}} 1 \frac{dxdy}{y^2} } \int_{ \Gamma_0(d) \backslash \mathbb{H}} |f(z) |^2 \frac{dxdy}{y^2}=1,
\end{align}
one has
\begin{align}
\label{rwc} \frac{1}{\int_{ \Gamma_0(d) \backslash \mathbb{H}} 1 \frac{dxdy}{y^2}} \int_{ \Gamma_0(d) \backslash \mathbb{H}} | f(z)|^4 \frac{dxdy}{y^2} \sim  \frac{1}{\sqrt{2\pi}}\intt t^{4} e^{\frac{-t^2}{2}} dt
\end{align}
as $\lambda_f \to \infty$ (and $d$ is fixed). The left hand side of (\ref{rwc}) is the fourth power of the $L^4$-norm of $f$, divided by the area of a fundamental domain. In the case $d=1$, Sarnak and Watson \cite[Theorem 3]{sar} have announced the upper bound $\lambda_f^\epsilon$ for the $L^4$-norm, possibly assuming the Ramanujan-Petersson conjecture, but the details have not yet appeared in print. Recently Luo \cite{luo} provided an elegant proof of the same upper bound in the case that $d$ is a prime discriminant and $f_\psi$ is the dihedral form ($L^2$-normalized) associated to a Grossencharacter $\psi$ of modulus 1 for $\mathbb{Q}(\sqrt{d})$. One would naturally be interested in going beyond this upper bound and obtaining an asymptotic for the fourth moment. By Parseval's theorem and spectral decomposition, we have that
\begin{align}
\label{specexp} \| f_\psi \|_4^4 = |\langle f_\psi^2, u_0 \rangle|^2 + \sum_{j\ge 1} |\langle f_\psi^2, u_j \rangle|^2 + (\text{continuous spectrum contribution}),
\end{align}
where $u_0$ is a constant and $\{ u_j : j\ge 1\}$ is an orthonormal Hecke-Maass basis for the cuspidal spectrum of $\Gamma_0(d)$ with trivial nebentypus. This consists of newforms of level 1 and of level d, and oldforms that are lifts from level 1. The continuous spectrum contribution is negligible (see \cite[section 5]{luo} and input the subconvexity bound from \cite[Theorem 1.1]{micven} instead of the convexity bound as Luo does). As for the cuspidal spectrum sum, as Luo explains in \cite[section 4]{luo}, the contribution of $u_j$ of level 1 or newforms of level d equals, by identities of Watson \cite{wat} and Ichino \cite{ich}, a mean of central values of $L$-functions having the shape
\begin{align}
\label{luomean}\sum_{0<t_j < T_{f_{\psi^2}}} \frac{1}{t_j {(T_{f_{\psi^2}}})^\half (1+T_{f_{\psi^2}}-t_j)^\half } \frac{L(\thalf, u_j) L(\thalf, u_j \times \chi_d) L(\thalf, u_j \times f_{\psi^2})}{L(1,f_{\psi^2})^2 L(1, \chi_d)^2 L(1,\sym^2 u_j)},
\end{align}
where $\chi_d$ is the real nebentypus of $f_\psi$, $\frac{1}{4}+t_j^2$ is the eigenvalue of $u_j$ and $\frac{1}{4} + (T_{f_{\psi^2}})^2$ is the eigenvalue of the form associated to the Grossencharacter $\psi^2$. 

The sum on the right hand side of (\ref{luomean}) may be divided into three parts: short ranges 
\begin{align}
\label{short2} 0<t_j< (T_{f_{\psi^2}})^{1-\epsilon}, 
\end{align}
and
\begin{align}
\label{short1} T_{f_{\psi^2}}-(T_{f_{\psi^2}})^{1-\epsilon}< t_j < T_{f_{\psi^2}}
\end{align}
on which the sum is expected to tend to 0 on the Lindel\"{o}f hypothesis, and the bulk range 
\begin{align}
\label{bulk} (T_{f_{\psi^2}})^{1-\epsilon} <t_j < T_{f_{\psi^2}}-(T_{f_{\psi^2}})^{1-\epsilon},
\end{align}
which is expected to yield the main term. For this reason, the bulk range may be the most interesting to study.

In this paper we prove an asymptotic with power saving for a mean value which is similar to (\ref{luomean}) in the bulk range (\ref{bulk}), where the normalization factor in the sum is of size about $(T_{f_{\psi^2}})^{-2}$. This mean value is a somewhat simplified version of (\ref{luomean}), meant as a `test' case on which to develop ideas that may be helpful in eventually proving an asymptotic for the fourth moment of dihedral forms.
\begin{theorem}\label{main} Let $\chi_d$ be a quadratic Dirichlet character of prime modulus $d\equiv 1 \bmod 4$ with $\chi_d(-1)=1$. Let $\{ u_j: j\ge 1\}$ denote an orthonormal basis of Hecke-Maass cusp forms for  $\text{SL}_2(\mathbb{Z})   \backslash  \mathbb{H}$ ordered by Laplacian eigenvalue $\frac{1}{4}+t_j^2$. Let $f$ be an even form from this basis with eigenvalue $\frac{1}{4}+T^2$, where $T>0$. There exists a computable $\delta > 0$ such that
\begin{align}
\label{mainline} \frac{1}{T^2} \sum_{j\ge 1} e^{\frac{-t_j^2}{T^2}}  \frac{L(\thalf, u_j) L(\thalf, u_j \times \chi_d)  L(\thalf, u_j \times f )}{ L(1,{\normalfont \text{sym}} ^2 u_j)} = \frac{2 L(1,\chi_d)L(1,f)L(1,f\times \chi_d)}{\pi^2 L(2,\chi_d)} \log T + C +O_d(T^{-\delta}),
\end{align}
where $C$ is a constant given in section \ref{diagsection}.
\end{theorem}
\noindent 
\noindent  For the problem of obtaining an asymptotic for the fourth moment of dihedral forms, our result must be worked out in greater generality. Firstly, we have used a simpler weight function $e^{-t_j^2/T^2}$ than what actually appears in the identities of Watson and Ichino. Secondly, we have taken all our forms $u_j$ and $f$ to be of level 1, while one must also consider $u_j$ of level $d$ (both newforms and oldforms lifted from level 1) and $f$ of fixed level $d$ and nebetypus $\chi_d$ (it is for this last difference that our main term is of size $\log T$ while the main term of (\ref{luomean}) should be a constant). Even after this, one has still to bound the sum (\ref{luomean}) over the short ranges, a problem which seems to to require methods different from those in this paper. We do not know how to treat the range (\ref{short2}), but the range (\ref{short1}) can be handled by applying H\"{o}lder's inequality as Luo does and then applying Jutila's \cite{jut} and Ivi\'c's \cite{ivi2} bounds for moments of $L(\half, u_j)$ in short intervals of $t_j$ close to $T$.

Apart from the connection to the $L^4$-norm problem, Theorem \ref{main} is interesting as a result in its own right. It offers an asymptotic with a power saving for a mixed moment (in the sense of \cite[Theorem 1.2]{bloetal}), in between Ivi\'c's \cite{ivi} asymptotic for the fourth moment of $L(\half, u_j)$ and the unestablished second moment of $L(\half, u_j \times f)$. The latter problem seems to be very difficult because were an asymptotic with a power saving known for it, one could presumably use an amplifier to obtain a subconvex bound for $L(\half, u_j\times f)$ for $|t_j-T|$ as small as $T^{1-\epsilon}$.  Thus Theorem \ref{main} seems to be at the edge of present methods. Indeed, its proof employs the full power of spectral theory and ultimately relies on a subconvex estimate for $GL(2)\times GL(2)$ $L$-functions.

\section{Sketch}

We give a very rough sketch to indicate the main ideas of the proof of Theorem \ref{main}. The notation is defined in the next section.

Using approximate functional equations, we write the left hand side of (\ref{mainline}) as 
\begin{align}
\label{sketch1} \frac{1}{T^2} \sum_{j\ge 1}  \frac{e^{\frac{-t_j^2}{T^2}}}{L(1,\sym^2 u_j)} \Big( \sum_{n< T^{1-\epsilon}} \frac{\lambda_j(n)}{n^\half} + \sum_{n< T^{1+\epsilon}} \frac{\lambda_j(-n)}{n^\half} \Big) \Big(2\sum_{m< T} \frac{\lambda_j(m)\chi_d(m)}{m^\half} \Big) \Big(2\sum_{r< T^2}\frac{\lambda_j(r) \lam(r)}{r^\half} \Big).
\end{align}
Here we use approximate functional equations for $ L(\thalf, u_j \times \chi_d)$ and $ L(\thalf, u_j \times f )$ which are valid for $u_j$ even, since otherwise $L(\thalf, u_j)$ vanishes. For $ L(\thalf, u_j)$ we use an uneven approximate functional equation. This is a useful idea which greatly simplifies the analysis.

The next step is to apply the Kuznetsov trace formula, and we must show that the off-diagonal part is bounded by a negative power of $T$. It is easily seen that the contribution of the shorter sum over $n$ in (\ref{sketch1}) is small, so that we are left to bound the other part of the off-diagonal (coming from Kuznetsov applied to opposite sign terms):
\begin{align}
\frac{1}{T^2} \sum_{\substack{n<T^{1+\epsilon}\\ m<T \\ r< T^2}}  \frac{\lam(r) \chi_d(m)}{\sqrt{nmr}} \sum_{c\ge 1} \frac{S(-nm, r, c)}{c}\intt \sinh(\pi t) K_{2it}\Big( \frac{\sqrt{nmr}}{c} \Big) e^{\frac{-t^2}{T^2}} tdt.
\end{align}
The Bessel transform is evaluated as a bump function of size $T$ supported on $\frac{\sqrt{nmr}}{c}\sim T$. Thus we need to bound
\begin{align}
\frac{1}{T^4} \sum_{\substack{c,n\sim T^{1+\epsilon}\\ m\sim T \\ r\sim  T^2}}  \lam(r) \chi_d(m)  S(-nm, r, c).
\end{align}
Ignoring the fixed character $\chi_d$ for the purposes of this sketch, Poisson summation in $n$ and $m$ (after splitting into residue classes modulo $c$) gives us roughly
\begin{align}
\frac{1}{T^2} \sum_{\substack{c\sim  T^{1+\epsilon}  \\ r\sim  T^2}} \  \sum_{a,b \bmod c} \ \sum_{|\ell_1|,|\ell_2|<T^{\epsilon}}  \frac{ \lam(r) S(-ab, r, c)}{c^2}e\Big(\frac{a\ell_1}{c}\Big) e\Big(\frac{b\ell_2}{c}\Big).
\end{align}
We consider the case $\ell_1=\ell_2=1$ and evaluate the exponential sum, getting
\begin{align}
\frac{1}{T^2} \sum_{\substack{c \sim  T^{1+\epsilon}  \\ r\sim T^2}} \  \summ_{x \bmod c} \  \frac{ \lam(r) }{c}e\Big(\frac{x(r+1)}{c}\Big).
\end{align}
Now Voronoi summation in $r$ gives us
\begin{align}
\label{sketch2}\frac{1}{T^2} \sum_{\substack{ q< T^2}}  \lam(q) \sum_{c\sim  T^{1+\epsilon}}    \frac{ S(q,1,c)}{c}.
\end{align}
We use Kuznetsov's formula to express the innermost sum of Kloosterman sums in terms of automorphic forms (we are in the nice range $c\asymp \sqrt{q}$). This reduces the proof to bounding by a negative power of $T$ the sum
\begin{align}
\frac{1}{T^2} \sum_{\substack{ q< T^2}}  \lam(q) \lambda_g(q) + \text{similar},
\end{align}
where $g$ is a Hecke-Maass cusp form with Laplacian eigenvalue bounded by $T^{\epsilon}$. The required estimate follows from a subconvex bound for $L(\half, f\times g)$.

\section{Preliminaries} \label{prelim}

\subsection*{Convention} Throughout the paper, $\epsilon>0$ denotes a small parameter which may be chosen to be as small as we like, but does not denote the same one from one occurrence to another. All implicit constants may depend on $\epsilon$ and $d$.

\subsection{\texorpdfstring{${\boldsymbol L}$-functions}{L-functions}}
 
 Let $\lambda_j(n)$ and $\lambda_f(n)$ denote the (real) eigenvalues of the $n$-th Hecke operator corresponding to $u_j$ and $f$ respectively, where we write $\lambda_j(-n)=\lambda_j(n)$ for $u_j$ even and $\lambda_j(-n)=-\lambda_j(n)$ for $u_j$ odd. The eigenvalues satisfy the multiplicative relations 
 \begin{align}
\label{hmult} \lambda_j(n)\lambda_j(m) = \sum_{l|(n,m)} \lambda_j\Big( \frac{nm}{l^2} \Big), \quad \lambda_j(mn) = \sum_{l | (n,m)} \mu(l) \lambda_j(\tfrac{m}{l}) \lambda_j(\tfrac{n}{l}),
 \end{align} 
where the sums above run over positive divisors only, the average bound
\begin{align}
\sum_{n\le x} |\lambda_j(n)| \ll x^{1+\epsilon},
\end{align}
and the individual bound $\lambda_j(n)\ll n^{\frac{7}{64}+\epsilon}$ of Kim-Sarnak \cite{kimsar}. At the infinite place, the Ramanujan-Petersson conjecture is known to be true: that is, $t_j$ is real. 

We have the $L$-functions
 \begin{align}
&L(s, u_j) = \sum_{n\ge 1} \frac{\lambda_f(n)}{n^s},\\
  &L(s, u_j \times \chi_d) = \sum_{n\ge 1} \frac{\lambda_f(n)\chi_d(n) }{n^s},
  \end{align}
 and
 \begin{align}
 \label{rsdef} L(s, u_j \times f) = \zeta(2s) \sum_{n\ge 1} \frac{\lambda_j(n) \lambda_f(n)}{n^s}
 \end{align}
 for $\Re(s)>1$ with analytic continuation to entire functions on the whole complex plane. The analytic conductors of the $L$-functions above are of size $1+|t_j|^2$, $1+|t_j|^2$, and $1+|T^2-t_j^2|^2$ respectively. Let $\Gam(s)= \pi^{-\frac{s}{2}} \Gamma(\frac{s}{2})$. For $u_j$ even we have the functional equations
\begin{align}
 \label{functeq1}  &L(s, u_j) \Gam(s+it_j)\Gam(s-it_j)= L(1-s, u_j) \Gam(1-s+it_j)\Gam(1-s-it_j),\\
 \label{functeq2}  &d^s L(s,u_j\times \chi_d) \Gam(s+it_j)\Gam(s-it_j)=  d^{1-s} L(1-s, u_j\times \chi_d) \Gam(1-s+it_j)\Gam(1-s-it_j),
\end{align}
and
\begin{multline}
 \label{functeq3} L(s, u_j\times f)  \Gam(s+it_j + iT)\Gam(s-it_j + iT)  \Gam(s+it_j- iT)\Gam(s-it_j- iT)  \\
= L(1-s, u_j\times f) \Gam(1-s+it_j + iT)\Gam(1-s-it_j + iT)  \Gam(1-s+it_j- iT)\Gam(1-s-it_j- iT)  .
\end{multline}
For $u_j$ odd we have the functional equation
\begin{align}
 \label{functeq4} L(s, u_j) \Gam(1+s+it_j)\Gam(1+s-it_j)= -L(1-s, u_j) \Gam(2-s+it_j)\Gam(2-s-it_j).
\end{align}
All of these may be found in \cite[chapters 3 and 7]{gol}.

Also define
\begin{align}
\label{eq:EisenHeckeDef} \lambda(n,t)=\sum_{ab=n} \Big(\frac{a}{b}\Big)^{it}.
\end{align}
These are the Hecke eigenvalues corresponding to the Eisenstein series. Note that $\lambda(n,t)$ satisfies the same Hecke relations (\ref{hmult}). For $\Re(s)>1$ we have
\begin{align}
\label{Etwist1}&\zeta(s-it)\zeta(s+it)=\sum_{n\ge 1}\frac{\lambda(n,t)}{n^s},\\
\label{Etwist2}&L(s-it,\chi_d)L(s+it,\chi_d) = \sum_{n\ge 1} \frac{\lambda(n,t)\chi_d(n)}{n^s},\\
\label{Etwist3}&L(s-it,f)L(s+it,f) = \zeta(2s) \sum_{n\ge 1} \frac{\lambda(n,t)\lambda_f(n)}{n^s}.
\end{align}
These identities can be seen by comparing Euler factors on both sides, as in \cite[section 3]{li}.

\subsection{Stirling's Approximation}\label{stirap}

For $t\gg 1$, $\sigma>0$ fixed and $0<\gamma <t^{\epsilon}$,  we have 
\begin{align}
\Gamma(\sigma +i\gamma + it)=\sqrt{2\pi}\exp\Big( (\sigma-\tfrac{1}{2}+i\gamma+it) \log (\sigma+i\gamma+it) - (\sigma+i\gamma+it)+ O(t^{-1+\epsilon})  \Big).
\end{align}
The complex logarithm equals
\begin{align}
 \log (\sigma+i\gamma+it) = \half\log \Big( (\gamma+t)^2+\sigma^2 \Big) + i \cot^{-1}\Big(\frac{\sigma}{\gamma+t} \Big)=\log  t  + \frac{\gamma}{t} +i\Big( \frac{\pi}{2} - \frac{\sigma}{t}\Big) + O(t^{-2+\epsilon}).
\end{align}
Thus 
\begin{align}
\Gamma(\sigma +i\gamma + it) = \sqrt{2\pi} t^{(\sigma-\half+i\gamma +it )} \exp\Big( i\tfrac{\pi}{2}(\sigma-\thalf+i\gamma +it)  -it + O(t^{-1+\epsilon})\Big).
\end{align}
Similarly,
\begin{align}
\Gamma(\sigma +i\gamma - it) = \sqrt{2\pi} t^{(\sigma-\half+i\gamma -it )} \exp\Big( - i\tfrac{\pi}{2}(\sigma-\thalf+i\gamma - it)  +it + O(t^{-1+\epsilon})\Big).
\end{align}
Of course this can be made more precise by taking more terms in Stirling's approximation and the Taylor series of $\log$ and $\cot^{-1}$ above.

\subsection{Approximate functional equations}

\begin{lemma} \label{afelemma} For any $\sigma>0$ and some parameter $0<\beta< \frac{1}{100}$ to be fixed later, let
\begin{align}
\label{eq:V1def}
V_1^{\pm}(x,t)=& \frac{1}{2\pi i} \int_{(\sigma)} e^{s^2} \left( x T^{\pm \beta} \right)^{-s}   \frac{ \Gam(\half +s + it) \Gam(\half  +s- it) }{ \Gam(\half + it)  \Gam(\half - it) }  \frac{ds}{s},\\
\label{eq:V0def}
V_1(x,t)=& \frac{1}{2\pi i} \int_{(\sigma)} e^{s^2}  x^{-s}   \frac{ \Gam(\half +s + it) \Gam(\half  +s- it) }{ \Gam(\half + it)  \Gam(\half - it) }  \frac{ds}{s},\\
\label{eq:V2def}
V_2 (x,t)=& \frac{1}{2\pi i} \int_{(\sigma)} e^{s^2} x^{-s}  d^s   \frac{ \Gam(\half +s + it) \Gam(\half  +s- it) }{ \Gam(\half + it)  \Gam(\half - it) }   \frac{ds}{s},\\
\label{eq:V3def}
V_3 (x,t)=& \frac{1}{2\pi i} \int_{(\sigma)} e^{s^2} \zeta(1+2s) x^{-s}  \prod_{\pm} \frac{ \Gam(\half +s + it \pm iT)\Gam(\half +s - it \pm iT)}{  \Gam(\half  + it \pm iT)  \Gam(\half - it \pm iT) }  \frac{ds}{s}.
\end{align}
For $T^{1-\epsilon}< t_j <T^{1+\epsilon}$, we have that
\begin{align}
\label{afe-error} L(\thalf, u_j)  = \sum_\pm \sum_{n\ge 1} \frac{\lambda_j(\pm n)}{n^{\half}}V_1^{\pm}(n,t_j) +O\big(T^{-\half+\frac{\beta}{2}+\epsilon}\big).
\end{align}
For $u_j$ even, we have that
\begin{align}
\label{afe-e1}&L(\thalf, u_j ) = 2\sum_{n\ge 1} \frac{\lambda_j(n) }{n^{\half}}V_1(n,t_j),\\
\label{afe-e2} &L(\thalf, u_j \times \chi_d) = 2\sum_{m\ge 1} \frac{\lambda_j(m) \chi_d(m)}{m^{\half}}V_2(m,t_j),\\ 
\label{afe-e3} &L(\thalf, u_j \times f) = 2\sum_{r\ge 1} \frac{\lambda_j(r)\lambda_f(r)}{r^\half}V_3(r,t_j).
\end{align}
\end{lemma}

\proof
These follow in a standard way from \cite[Theorem 5.3]{iwakow} (by putting $G(u)=e^{u^2}$ and $X=1$ there) and the functional equations (\ref{functeq1}-\ref{functeq3}), but (\ref{afe-error}) requires some explanation. We start with the approximate functional equation (which follows from \cite[Theorem 5.3]{iwakow} by putting $G(u)=e^{u^2}$ and $X=T^\beta$ there),
\begin{align}
\label{uba} L(\thalf, u_j)  = \sum_\pm\sum_{n\ge 1} \frac{\lambda_j( \pm n)}{n^{\half}}  \frac{1}{2\pi i} \int_{(\sigma)} e^{s^2} \left( x T^{\pm \beta} \right)^{-s}   \frac{ \Gam(\half +s +\kappa_j+ it) \Gam(\half  +s+ \kappa_j- it) }{ \Gam(\half + \kappa_j+ it)  \Gam(\half + \kappa_j- it) }  \frac{ds}{s},
\end{align}
where $\kappa_j=0$ or $1$ as $u_j$ is even or odd. By the rapid decay of $e^{{s^2}}$ in vertical lines, we may restrict the integral above to $|\Im(s)|<T^{\epsilon}$. By Stirling's approximation, for $\Re(s)>0$ fixed, $|\Im(s)|<T^{\epsilon}$ and $T^{1-\epsilon}< t <T^{1+\epsilon}$, we have
\begin{align}
 \frac{ \Gam(\half +s +1+ it) \Gam(\half  +s+ 1- it) }{ \Gam(\half + 1+ it)  \Gam(\half + 1- it) } =  \frac{ \Gam(\half +s+ it) \Gam(\half  +s- it) }{ \Gam(\half + it)  \Gam(\half - it) } + O(T^{-1+\epsilon})
\end{align}
 Thus up to a small error, the ratio of Gamma functions in (\ref{uba}) does not depend on $\kappa_j$, and (\ref{afe-error}) follows.
\endproof

We describe a trick that we will use. By (\ref{afe-e2}) and (\ref{afe-e3}), we have that 
\begin{multline}
\label{trick} L(\thalf, u_j)L(\thalf, u_j \times \chi_d)L(\thalf, u_j \times f) \\= L(\thalf, u_j)\left(2\sum_{m\ge 1} \frac{\lambda_j(m) \chi_d(m)}{m^{\half}}V_2(m,t_j)\right)\left( 2\sum_{r\ge 1} \frac{\lambda_j(r)\lambda_f(r)}{r^\half}V_3(r,t_j)\right)
\end{multline}
holds for even forms $u_j$. But when $u_j$ is odd, $L(\thalf, u_j )=0$ and both sides vanish. So the equality holds for odd forms too. For the factor $L(\thalf, u_j )$ on the right hand side, we may use the uneven approximate functional equation given by (\ref{afe-error}).

On the critical line, we will need the following approximate functional equations.
\begin{lemma} \label{afelem2} Keeping with the notation of Lemma \ref{afelemma}, we have
\begin{align}
&|\zeta(\thalf +it)|^2  = \sum_\pm \sum_{n\ge 1} \frac{\lambda(n,t)}{n^{\half}}V_1^{\pm}(n,t),\\
&|L(\thalf + it,\chi_d)|^2 = 2\sum_{m\ge 1} \frac{\lambda(m,t)\chi_d(m)}{n^{\half}}V_2(m,t),\\
&|L(\thalf+it,f)|^2 = 2\sum_{m\ge 1} \frac{\lambda(r,t) \lambda_f(r) }{n^{\half}}V_3(r,t).
\end{align}
\proof
These follow by (\ref{Etwist1}-\ref{Etwist3}), \cite[Theorem 5.3]{iwakow}, and the functional equations of the relevant $L$-functions.
\endproof

\end{lemma}

\subsection{Kuznetsov trace formula}
We define the Kloosterman sums
\begin{align}
S(n,m,c) = \summ_{a \bmod c} e\left(\frac{na+m\bar{a}}{c}\right)
\end{align}
and, for $d|c$,
\begin{align}
S_{\chi_d}(n,m,c) = \summ_{a \bmod c} \chi_d(a) e\left(\frac{na+m\bar{a}}{c}\right)
\end{align}
where $e(x) =e^{2\pi i x}$, the sum is restricted primitive residue classes and $\bar{a}a \equiv 1 \bmod{c}$.

We recall Kuznetsov's trace formula. The spectral side of this formula is usually written in terms of the Fourier coefficients of an orthonormal basis of cusp forms, but we write it in terms of Hecke eigenvalues, using the relationship
\begin{align}
&u_j(x+iy)= \rho_j(1)\sum_{n\neq 0} \lambda_j(n) \sqrt{y} K_{i t_j} (2\pi n y) e(nx),\\
\label{rho1formula} &\rho_j(1)^2 = \frac{2 \cosh(\pi t_j)}{L(1,\sym^2 u_j)}.
\end{align}
The calculation for (\ref{rho1formula}) may be found in \cite[section 3]{blo}. Let
\begin{align}
\mathcal{J}^+(x, t) = \frac{2i}{\sinh(\pi t)} J_{2it}(4\pi x), \quad \mathcal{J}^-(x, t) = \frac{4}{\pi} K_{2it}(4\pi x) \cosh(\pi t).
\end{align}
We have
\begin{lemma} \label{lem:Kuznetsov} \cite[Theorems 2.2, 2.4]{moto}
Let $h(z)$ be an even, holomorphic function on $|\Im(z)| < \frac{1}{2}+\theta$ with decay $|h(z)| \ll (1+|z|)^{-2-\theta}$ on that strip, for some $\theta>0$. Then for $n,m > 0$,
\begin{multline}
\sum_{j \ge 1} \frac{\lambda_j(\pm n) \lambda_j(m)}{L(1,{\normalfont \text{sym}} ^2 u_j)} h(t_j) +\int_{-\infty}^\infty \frac{\lambda(n, t) \lambda(m, -t)}{|\zeta(1+2it)|^2} h(t) \frac{dt}{2\pi} \\
= \delta_{\pm n,m} \int_{-\infty}^\infty h(t) \frac{d^*t}{2\pi^2} + \sum_{c \ge 1} \frac{S(\pm n,m, c)}{c} \int_{-\infty}^\infty \mathcal{J}^\pm(\tfrac{\sqrt{nm}}{c}, t) h(t) \frac{d^* t}{2 \pi},
\end{multline}
where $\delta_{n,m}$ is 1 if $n=m$ and 0 otherwise, $\delta_{-n,m}$ is always 0, and $d^* t =  \tanh(\pi t) \, t dt$.
\end{lemma}

We will also need Kuznetsov's formula from the geometric side to the spectral side, written in terms of Hecke eigenvalues again. Let $\mathcal{B}_k(d,\chi_d)$ denote the orthonormal basis of Hecke eigenforms for the space of holomorphic cusp forms of weight $k$ for $\Gamma_0(d)$ with nebentypus $\chi_d$. Let $\mathcal{B}(d,\chi_d)$ denote the orthonormal basis of Hecke eigenforms for the space of Maass cusp forms for $\Gamma_0(d)$ with nebentypus $\chi_d$. Recall from \cite[page 373]{iwakow} that $\mathcal{B}(d,\chi_d)$ and $\mathcal{B}_k(d,\chi_d)$ consist of newforms, since $\chi_d$ is primitive. For $g\in \mathcal{B}_k(d,\chi_d)$ or $g\in\mathcal{B}(d,\chi_d)$, let $\lambda_g(n)$ denote the eigenvalue of the $n$-th Hecke operator corresponding to $g$. For $g\in \mathcal{B}(d,\chi_d)$, let $\half+t_g^2$ denote its Laplacian eigenvalue. Let $\rho_g(1)$ be first Fourier coefficient of $g$ as defined in \cite[Section 2.1.3]{bloharmic}, a normalization factor to go between the Fourier coefficients and the Hecke eigevalues $\lambda_g(n)$. We will need the bounds (see \cite[Section 2.6]{harmic}; note the slightly different definition of $\rho_g(1)$ given there for holomorphic forms):
\begin{align}
\label{rho1h} \frac{(4\pi)^{k-1}}{(k-1)!}k^{-\epsilon} \ll |\rho_g(1)|^2\ll \frac{(4\pi)^{k-1}}{(k-1)!}k^\epsilon
\end{align}
for $g\in  \mathcal{B}_k(d,\chi_d)$ and
\begin{align}
\label{rho1m} \cosh(\pi t_g) {(1+|t_g|)}^{-\epsilon} \ll |\rho_g(1)|^2\ll \cosh(\pi t_g) {(1+|t_g|)}^{\epsilon} 
\end{align}
for $g\in \mathcal{B}(d,\chi_d)$.

Let $E_{v}(z,s)$ denote the Eisenstein series associated with the singular cusp $\frac{1}{v}$, for $v|d$. Its $n$-th Fourier coefficient can be written in terms of 
\begin{align}
\label{lameis} \lambda_v(n,t) = \sum_{ab=n} \chi_v(a) \chi_{\frac{d}{v}}(b) \Big(\frac{a}{b}\Big)^{it}, 
\end{align}
where $\chi_v\chi_{\frac{d}{v}} = \chi_d$, with a go between factor $\rho_v(1)$ that satisfies the same bound as above: 
\begin{align}
\cosh(\pi t) {(1+|t|)}^{-\epsilon} \ll |\rho_v(1)|^2\ll \cosh(\pi t) {(1+|t|)}^{\epsilon}.
\end{align}
These facts can be found in \cite[sections 6 and 7]{dfi}.

\begin{lemma} \label{kuzbackwards} \cite[Section 2.1.4]{bloharmic} Let $\Phi$ be a smooth function compactly supported on the positive real numbers. Let
\begin{align}
&\dot{\Phi}(k) = i^k\int_0^\infty J_{k-1}(x)\Phi(x) \frac{dx}{x},\\
\label{phitilde} &\tilde{\Phi}(t)= \frac{i}{2\sinh(\pi t)} \int_0^\infty (J_{2it}(x) - J_{-2it}(x))\Phi(x)\frac{dx}{x},\\
&\check{\Phi}(t)= \frac{2\cosh(\pi t)}{\pi} \int_0^\infty K_{2it}(x) \Phi(x) \frac{dx}{x}.
\end{align}
For positive integers $q$ and $\ell$, we have
\begin{align}
\label{kuzbackpos} \sum_{c\ge 1} \frac{S_{\chi_d}(q,\ell,cd)}{cd} \Phi\left( \frac{4\pi \sqrt{q \ell }}{cd}  \right) &= \sum_{\substack{k > 0\\ k\equiv 0 \bmod 2}} \sum_{g\in \mathcal{B}_k(d,\chi_d) } \dot{\Phi}(k) \frac{(k-1)!}{\pi (4\pi)^{k-1}} |\rho_g(1)|^2 \lambda_g(q)\lambda_g(\ell)\\
\nonumber &+\sum_{g\in \mathcal{B}(d, \chi_d)} \tilde{\Phi}(t_g) \frac{4\pi}{\cosh(\pi t_g)} |\rho_g(1)|^2 \lambda_g(q)\lambda_g(\ell) \\
\nonumber &+\sum_{v|d} \intt \tilde{\Phi}(t) \frac{4\pi}{\cosh(\pi t)} |\rho_{v}(1)|^2 \lambda_{v}(q,t)\lambda_{v}(\ell,-t) dt
\end{align}
and
\begin{align}
\label{kuzbackneg} \sum_{c\ge 1} \frac{S_{\chi_d}(q,-\ell,cd)}{cd} \Phi\left( \frac{4\pi \sqrt{q \ell }}{cd}  \right) &=\sum_{g\in \mathcal{B}(d, \chi_d)} \check{\Phi}(t_g) \frac{4\pi}{\cosh(\pi t_g)} |\rho_g(1)|^2 \lambda_g(q)\lambda_g(\ell) \\
\nonumber &+\sum_{v|d} \intt \check{\Phi}(t) \frac{4\pi}{\cosh(\pi t)} |\rho_{v}(1)|^2 \lambda_{v}(q,t)\lambda_{v}(\ell,-t) dt.
\end{align}
\end{lemma}
\noindent Suppose that $\Phi$ is a smooth function compactly supported on $(T^{-\epsilon},T^{\epsilon})$, with derivatives satisfying $\| \Phi^{(n)} \|_\infty \ll (T^\epsilon)^n$. Then we note that the sums in (\ref{kuzbackpos}) may effectively be restricted to $k<T^{\epsilon}$, $|t_g|< T^{\epsilon}$ and $|t|<T^\epsilon$, as the contribution of the larger parameters is less than $T^{-100}$, say. For the holomorphic forms, this may be seen by the bound \cite[8.402]{gr}
\begin{align}
J_{k-1}(x) \ll \frac{1}{\Gamma(k)} \Big(\frac{x}{2} \Big)^{k-1},
\end{align}
valid for $x\in (T^{-\epsilon}, T^{\epsilon})$ and $k>T^{3\epsilon}$. For the Maass forms and Eisenstein series, this may be seen by repeatedly integrating by parts the integral in (\ref{phitilde}) for $x\in (T^{-\epsilon},T^{\epsilon})$ and $t >T^{3\epsilon}$, after applying the power series expansion  \cite[8.402]{gr}
\begin{align}
\frac{J_{2it}(x)}{2\sinh (\pi t)} = \frac{1}{2\sinh (\pi t)} \sum_{n\ge 0} \frac{(-1)^n x^{2n+2it}}{n!\Gamma(n+2it+1)}
\end{align}
which converges absolutely. Similarly, the sums in (\ref{kuzbackneg}) may be restricted to $|t_g|<T^{\epsilon}$ and $|t|<T^{\epsilon}$.

When $|k_g|,|t_g|, |t| <T^{\epsilon}$ and $T^{-\epsilon}< x < T^{\epsilon}$, we have that
\begin{align}
J_{k-1}(x), \  \frac{J_{2it}(x)}{\sinh (\pi t)} , \ \cosh(\pi t) K_{2it}(x) \ll T^{\epsilon}
\end{align}
by \cite[8.411 1]{gr}, \cite[lemma 6]{blomil2} and \cite[proposition 9]{harmic}.

\subsection{Subconvexity}

We record the following subconvex bounds, for some $\Delta>0$:
\begin{align}
L(\thalf + it, f \times g) \ll T^{1-\Delta},
\end{align}
where $|t|<T^{\epsilon}$, $|t_g|<T^{\epsilon}$ if $g\in \mathcal{B}(d,\chi_d)$ and $k<T^{\epsilon}$ if $g\in \mathcal{B}_k(d,\chi_d)$, and
\begin{align}
\label{eq:fSubconv}
L(\thalf + it, f), \quad L(\thalf + it, f \times \chi_d) \ll T^{\half-\Delta},
\end{align}
where $|t|<T^{\epsilon}$. These bounds may be found in \cite[Theorems 1.1 and 1.2]{micven}, which provide sufficiently general results that allow nontrivial nebentypus and bounds which depend polynomially on $|t_g|$, $k$, and $|t|$. It follows in a standard way, using Perron's formula, that for some $\Delta>0$ we have
\begin{align}
\label{subconvsum1} &\sum_{n<N} \lambda_f(n)\lambda_g(n) \ll N^\half T^{1-\Delta},\\
\label{subconvsum2} &\sum_{n<N} \lambda_f(n)\chi_d(n)  \ll N^\half T^{\half-\Delta}.
\end{align}
Note that for $\Re(s)>1$ and $v|d$, we have
\begin{align}
L(s+it,f\times \chi_v)L(s - it, f\times \chi_{\frac{v}{d}}) = \sum_{n\ge 1} \frac{\lambda_f(n)\lambda_v(n,t)}{n^s}.
\end{align}
This follows by comparing Euler products on both sides, as in (\ref{Etwist1}-\ref{Etwist3}). Thus 
\begin{align}
\label{subconvsum3} \sum_{n<N} \lambda_f(n)\lambda_v(n,t) \ll N^\half T^{1-\Delta}
\end{align}
for some $\Delta>0$.

\subsection{Spectral large sieve}

\begin{lemma} \cite{jut2} \label{largesieve}  For an arbitrary complex sequence $\{ a_n \}$, we have
\begin{align}
\sum_{|t_j-T|<A} \Big| \sum_{n< N} a_n \lambda_j(n) \Big|^2 \ll (TN)^{\epsilon}(AT+N)\Big(\sum_{n<N} |a_n|^2\Big).
\end{align} 
\end{lemma}

\subsection{Voronoi summation formula}

\begin{lemma} \label{voronoilem} \cite[Theorem 4.2]{godber}
Let $\phi$ be a smooth function with compact support on $(1,2)$. For $(a,c)=1$, we have 
\begin{align}
\sum_{r\ge 1}\lambda_f(r)e\left(\frac{r\overline{a}}{c}\right) \phi\left(\frac{r}{R}\right) =c \sum_\pm \displaystyle\sum_{q \ge 1}  \frac{\lambda_f(q)}{q}e\left(\frac{\pm qa}{c}\right)  \int_{(\sigma)} \tilde{\phi}(-s)  \left(\frac{\pi^2 Rq}{c^2}\right)^{-s} G_\pm(s) ds,
\end{align}
where $\sigma>-1$, $\tilde{\phi}$ is the Mellin transform of $\phi$ and 
\begin{align}
 G_\pm(s) 4\pi^2 i = \frac{ \Gamma\left(\frac{1+s+iT}{2}\right) \Gamma\left(\frac{1+s-iT}{2}\right) }{ \Gamma\left(\frac{-s+iT}{2}\right) \Gamma\left(\frac{-s-iT}{2}\right)} \pm \frac{ \Gamma\left(\frac{2+s+iT }{2}\right) \Gamma\left(\frac{2+s-iT }{2}\right) }{ \Gamma\left(\frac{1-s+iT}{2}\right) \Gamma\left(\frac{1-s-iT }{2}\right)}.
\end{align}
\end{lemma}

\subsection{Averages of Bessel functions}

\begin{lemma} \label{jbeslemma}
We have that
\begin{align}
\label{jbeslem} \intt \frac{J_{2it}(2\pi x)}{\cosh(\pi t)} h\Big(\frac{t}{T}\Big) t dt = &\frac{-i\sqrt{2}}{\pi} \frac{T^2}{\sqrt{x}} \Re\left((1+i)e(x)   \int_0^\infty t h(t) e\Big(\frac{-t^2 T^2}{2\pi^2 x}\Big) dt \right) \\
\nonumber &+ O\Big( \frac{x}{T^{3-12\alpha}} \Big) +O(T^{-100})
\end{align}
for any $x>0$ and any smooth even  function $h$, compactly supported on $(T^{-\alpha}, T^{\alpha})\cup (-T^{-\alpha}, -T^{\alpha})$ with derivatives satisfying $\| h^{(n)} \|_\infty \ll (T^{\alpha})^{n}$ for some $0<\alpha< \frac{1}{100}$. The main term is bounded by $T^{-100}$ if $x<T^{2-3\alpha}$.
\end{lemma}
\proof
We follow the ideas in \cite[Lemma 5.8]{iwaniec}. By \cite[pg. 180]{watson}, we have that
\begin{align}
 \frac{J_{2it}(2\pi x) - J_{-2it}(2\pi x)}{\cosh \pi t} = -2i \tanh (\pi t) \intt  \cos(2\pi x\cosh\pi u - 2\pi tu) du.
\end{align}
So the left hand side of (\ref{jbeslem}) equals
\begin{align}
\label{pois11}  -2 i \Re\left( \intt  e(x\cosh \pi u )   \int_0^\infty    \tanh (\pi t)  h\Big(\frac{t}{T}\Big) t  e(-ut) dt \ du \right).
\end{align}
We may replace $ \tanh \pi t$ by 1, with an admissible error since $ \tanh \pi t=1+O(e^{-t})$. Then by integrating by parts several times the $t$-integral, we see that the contribution of $|u|>T^{-1+2\alpha}$ is less than $T^{-100}$, say. For $|u|\le T^{-1+2\alpha}$, we take the Taylor expansion of $\cosh \pi u$. Following these steps, we see that (\ref{pois11}) equals
\begin{align}
\label{pois12} -2i \Re\left( e(x) \intt e\Big(\frac{x(\pi u)^2}{2} \Big)  \int_0^\infty   h\Big(\frac{t}{T}\Big) t  e(-ut) dt  \ du \right) + O\Big( \frac{x}{T^{3-12\alpha}} \Big).
\end{align}
Now using that
\begin{align}
 \intt e(u^2 y) e(- t u) du = \frac{1+i}{2\sqrt{y}}e\Big(\frac{-t^2}{4y}\Big),
\end{align}
we have that (\ref{pois12}) equals
\begin{align}
\frac{-i\sqrt{2}}{\pi} \frac{T}{\sqrt{x}} \Re\left((1+i)e(x)   \int_0^\infty \frac{t}{T}  h\Big(\frac{t}{T}\Big) e\Big(\frac{-t^2}{2\pi^2 x}\Big) dt \right).
\end{align}
Repeated integration by parts shows that the integral above is less than $T^{-100}$, say, if $x<T^{2-3\alpha}$.
\endproof

\begin{lemma}\label{kbeslemma} We have that
\begin{align}
\label{kbeslem}  \intt  \sinh(\pi t) K_{2it}(2\pi x) h\Big(\frac{t}{T}\Big) t dt =  \frac{\pi T}{2} H\Big( \frac{\pi x}{T} \Big) - \frac{i\pi^3}{12T}H^{(3)}\Big( \frac{\pi x}{T} \Big)+ O\Big( \frac{x}{T^{4-14\alpha}} \Big) + O(T^{-100}),
\end{align}
where $H(y)=yh(y)$, for any $x>0$ and any smooth even  function $h$, compactly supported on $(T^{-\alpha}, T^{\alpha})\cup (-T^{-\alpha}, -T^{\alpha})$ with derivatives satisfying $\| h^{(n)} \|_\infty \ll (T^{\alpha})^{n}$ for some $0<\alpha< \frac{1}{100}$.
\end{lemma}
 
\proof
Again, we basically follow the ideas in \cite[Lemma 5.8]{iwaniec}.
By \cite[8.432 4]{gr} we have that
\begin{align}
\sinh(\pi t) K_{2it}(2\pi x) = \frac{\pi \tanh (\pi t)}{2} \intt \cos (2\pi x \sinh \pi u) e(tu) du.
\end{align}
So the left hand side of (\ref{kbeslem}) equals
\begin{align}
\label{k01}
\frac{\pi}{2} \intt \cos (2\pi x \sinh \pi u) \intt e(tu) \tanh (\pi t) h\Big(\frac{t}{T}\Big) t dt du.
\end{align}
The inner $t$-integral is even function of $u$. Therefore (\ref{k01}) equals
\begin{align}
\label{k02} \frac{\pi}{2} \intt e( -x \sinh \pi u) \intt e(tu) \tanh (\pi t) h\Big(\frac{t}{T}\Big) t dt du.
 \end{align}
Integrating by parts several times the $t$-integral shows that the contribution of $|u|>T^{-1+2\alpha}$ is less than $T^{-100}$. For $|u|\le T^{-1+2\alpha}$, we take the Taylor expansion of $\sinh \pi u$, getting that (\ref{k02}) equals
\begin{multline}
\frac{\pi}{2} \intt e(-\pi x u)  \intt e(tu)  \tanh (\pi t) h\Big(\frac{t}{T}\Big) t dt du \\
+  \frac{\pi^3 x}{12} \intt  u^3 e(- \pi x u)  \intt e(tu)  \tanh (\pi t)  h\Big(\frac{t}{T}\Big) t dt du  + O\Big( \frac{x}{T^{4-14\alpha}} \Big).
\end{multline}
We may replace $ \tanh \pi t$ by 1, with an admissible error since $ \tanh \pi t=1+O(e^{-t})$. Then by Fourier inversion, the main term equals
\begin{align}
 \frac{\pi^2 x}{2}   h\Big(\frac{\pi x}{T}\Big) - \frac{i \pi^2 x}{12}\frac{d^3}{dx^3}\Big( x  h\Big(\frac{\pi x}{T}\Big)\Big).
\end{align}
 \endproof

\subsection{Test functions}

Define, for $0<\alpha<\frac{1}{100}$,
\begin{align}
&W_1(t)=\exp\left(\frac{-t^2}{T^2}\right),\\
\label{defw2} &W_2(t) = \left(1-\exp \left(-\left(\frac{t}{T^{1-\frac{\alpha}{2}}}\right)^{2\lceil \frac{1000}{\alpha} \rceil} \right) \right) \left(1-\exp \left(-\left(\frac{T^2-t^2}{T^{2-\frac{\alpha}{2}}}\right)^{2\lceil \frac{1000}{\alpha} \rceil} \right) \right).
\end{align}
By taking $\alpha$ small enough, we have that $W_1(t)W_2(t)$ is less than $T^{-100}$ unless 
\begin{align}
\label{wrange}  \Big| \frac{t}{T} \Big| \in (T^{-\alpha} ,1 - T^{-\alpha} ) \cup  (1 + T^{-\alpha}, T^{\epsilon}  ),
\end{align}
in which range 
\begin{align}
T^n \frac{d^n}{dt^n} W_1(t)W_2(t) \ll \left( T^{\alpha} \right)^n.
\end{align}
The point of these weight functions is that they are designed to satisfy the conditions of Kuznetsov's trace formula and to localize $t$ near $T$, but not too near so as to cause conductor-dropping of the Rankin-Selberg $L$-function (\ref{rsdef}). We have that the left hand side of (\ref{mainline}) equals
\begin{multline}
\frac{1}{T^2} \sum_{j\ge 1}  W_1(t_j) W_2(t_j) \frac{L(\thalf, u_j) L(\thalf, u_j \times \chi_d)  L(\thalf, u_j \times f )}{ L(1,{\normalfont \text{sym}} ^2 u_j)} \\+ \frac{1}{T^2} \sum_{j\ge 1}  W_1(t_j) (1-W_2(t_j)) \frac{L(\thalf, u_j) L(\thalf, u_j \times \chi_d)  L(\thalf, u_j \times f )}{ L(1,{\normalfont \text{sym}} ^2 u_j)}.
\end{multline}
The second sum above may be restricted to $u_j$ even, since $L(\thalf,u_j)=0$ otherwise, and to $\mathcal{S}=\{ t_j <T^{1-\alpha} \}\cup \{T-T^{1-\alpha} < t_j <T+T^{1-\alpha} \} $ by definition (\ref{defw2}). So by H\"{o}lder's inequality and the bound $L(1,{\normalfont \text{sym}} ^2 u_j)\gg T^{-\epsilon}$, we have that this sum is less than
\begin{align}
\label{holder} T^\epsilon \left( \frac{1}{T^2} \sum_{\substack{t_j \in \mathcal{S} \\ u_j \text{ even }}}  |L(\thalf, u_j)|^4   \right)^\frac{1}{4} \left( \frac{1}{T^2} \sum_{\substack{t_j \in \mathcal{S} \\ u_j \text{ even }}} |L(\thalf, u_j \times \chi_d)|^4  \right)^\frac{1}{4}\left( \frac{1}{T^2}\sum_{\substack{t_j \in \mathcal{S} \\ u_j \text{ even }}}  |L(\thalf, u_j \times f )|^2   \right)^\frac{1}{2} .
\end{align}
By (\ref{afe-e1}-\ref{afe-e3}), we see that the series for $|L(\thalf, u_j)|^2$ and $|L(\thalf, u_j\times \chi)|^2$ have length at most $t_j^{2+\epsilon}$, and the series for $L(\thalf, u_j \times f )$ has length at most $T^{\epsilon}(1+|t_j^2-T^2|)$. Thus by Lemma \ref{largesieve} we have that (\ref{holder}) is less than a negative power of $T$. It suffices therefore to give an asymptotic for
\begin{align}
\label{sumA} \frac{1}{T^2} \sum_{j\ge 1}  W_1(t_j) W_2(t_j) \frac{L(\thalf, u_j) L(\thalf, u_j \times \chi_d)  L(\thalf, u_j \times f )}{ L(1,{\normalfont \text{sym}} ^2 u_j)}. 
\end{align}

 In (\ref{eq:V1def}-\ref{eq:V3def}), write $s=\sigma+i\gamma$ for $\sigma>0$ and note that the integrals may be restricted to $|\gamma|<T^\epsilon$ by the rapid decay of $e^{s^2}$ in vertical lines. We restrict to this range of $\gamma$ and the range (\ref{wrange}) of $t$. By Stirling's approximation (see section \ref{stirap} for more details), we have 
\begin{align}
\label{v1approx} &V_1^\pm(x,t)= \frac{1}{2\pi i} \int_{(\sigma)} e^{s^2} \left( \frac{2\pi  x T^{\pm \beta}  }{|t|} \right)^{-s} \frac{ds}{s} + O(T^{-1+\alpha+\epsilon}) ,\\
\label{v2approx} &V_2(x,t)= \frac{1}{2\pi i} \int_{(\sigma)} e^{s^2} \left( \frac{2\pi  d x  }{|t|} \right)^{-s}   \frac{ds}{s} + O(T^{-1+\alpha+\epsilon}),
\end{align}
and
\begin{align}
\label{v3approx} V_3(x,t) &= \frac{1}{2\pi i} \int_{(\sigma)} e^{s^2} \zeta(1+2s) \left( \frac{4\pi^2 x  }{|T^2-t^2|} \right)^{-s} \frac{ds}{s} + O(T^{-1+\alpha+\epsilon}).
\end{align}
This can be made more precise by taking more terms in Stirling's approximation. We have
\begin{align}
\label{v1approx0} &V_1^\pm(x,t)= \frac{1}{2\pi i} \int_{(\sigma)} e^{s^2} \left( \frac{2\pi x T^{\pm \beta}   }{|t|} \right)^{-s} \Big( 1 + \sum_{n=1}^{1000} \frac{C_n(\sigma,\gamma)}{|t|^{n}} \Big)\frac{ds}{s} +O(T^{-100}),\\
\label{v2approx0} &V_2(x,t)= \frac{1}{2\pi i} \int_{(\sigma)} e^{s^2} \left( \frac{2\pi  d x }{|t|} \right)^{-s} \Big( 1 + \sum_{n=1}^{1000} \frac{C_n(\sigma,\gamma)}{|t|^{n}} \Big)  \frac{ds}{s} + O(T^{-100}),
\end{align}
and
\begin{multline}
\label{v3approx0} V_3(x,t) = \frac{1}{2\pi i} \int_{(\sigma)} e^{s^2} \zeta(1+2s) \left( \frac{4\pi^2  x }{|T^2-t^2|} \right)^{-s}\Big( 1 + \sum_{n=1}^{1000} \frac{C_n(\sigma,\gamma)}{|T-t|^{n}} \Big)\Big( 1 + \sum_{n=1}^{1000} \frac{C_n(\sigma,\gamma)}{|T+t|^{n}} \Big) \frac{ds}{s} \\
+ O(T^{-100}),
\end{multline}
for some $C_n(\sigma,\gamma)$ (not necessarily the same in each expression above) polynomial in $\sigma$ and $\gamma$.

Define
\begin{multline}
V^{\pm}\left(x_1,x_2,x_3; y \right) \\=  \int_{(\sigma)} e^{s_1^2} ( 2\pi x_1 T^{\pm \beta} )^{-s_1} |y|^{s_1}   \frac{ds_1}{s_1}  \cdot \int_{(\sigma)} e^{s_2^2} (2\pi d x_2 ) ^{-s_2} |y|^{s_2}   \frac{ds_2}{s_2} \cdot \int_{(\sigma)} e^{s_3^2} (4\pi^2 x_3  )^{-s_3} \zeta(1+2s_3) |1-y^2|^{-s_3}   \frac{ds_3}{s_3},
\end{multline}
for any $\sigma>0$.
Let $Z$ be any smooth, even function compactly supported on 
\begin{align}
\label{zsupp} (-T^{\epsilon}, -1 - T^{-\alpha}  )\cup (-1 + T^{-\alpha} , -T^{-\alpha} ) \cup   (T^{-\alpha} ,1 - T^{-\alpha} ) \cup  (1 + T^{-\alpha}, T^{\epsilon}  )
\end{align}
 with derivatives satisfying
\begin{align}
\| Z^{(n)}\|_\infty \ll \left( T^{\alpha} \right)^n.
\end{align}

\section{Applying the trace formula}

By Lemma \ref{afelemma} and (\ref{trick}), we have that (\ref{sumA}) equals
\begin{align}
\label{applyingtrace} & \frac{4}{T^2} \sum_{\pm}  \sum_{j\ge 1}  W_1(t_j) W_2(t_j) \sum_{n,m,r\ge 1} \frac{ \lambda_j(\pm n)\lambda_j(m)\chi_d(m)\lambda_j(r)\lambda_f(r)}{ L(1,{\normalfont \text{sym}} ^2 u_j) \sqrt{nmr}} V_1^\pm(n,t_j)V_2(m,t_j)V_3(r,t_j)\\
\nonumber &+O\Bigg(\frac{1}{T^{\frac{5}{2}-\frac{\beta}{2}-\epsilon} } \sum_{t_j<T^{1+\epsilon}} |L(\thalf, u_j \times \chi_d)  L(\thalf, u_j \times f ) | \Bigg).
\end{align}
The error term here arises from the error term of (\ref{afe-error}). It can easily be bounded by $T^{-\frac{1}{2}+\frac{\beta}{2}+\epsilon}$ on using Cauchy-Schwarz and Lemmas \ref{afelemma} and \ref{largesieve} (which amounts to the Lindel\"{o}f bound on average).

Using (\ref{hmult}), we have that the main term of (\ref{applyingtrace}) equals
\begin{align}
\label{eq:PreKuz}
 \frac{4}{T^2} \sum_{\pm} \sum_{j\ge 1}  W_1(t_j) W_2(t_j) \sum_{k,n,m,r\ge 1} \frac{\lambda_j(\pm n m)\chi_d(k m)\lambda_j(r)\lambda_f(r)}{ L(1,{\normalfont \text{sym}} ^2 u_j) k \sqrt{nmr}} V_1^\pm(k n,t_j)V_2(k m,t_j)V_3(r,t_j).
\end{align}
Now applying Lemma \ref{lem:Kuznetsov}, for each sign $\pm$, with
\begin{align}
 h(t) = W_1(t) W_2(t) V_1^\pm(k n,t)V_2(k m,t)V_3(r,t), 
\end{align}
we have that \eqref{eq:PreKuz} can be divided into three parts: The diagonal part with $r=nm$,
\begin{align}
\label{eq:KuzDiag}
 \frac{4}{T^2}  \sum_{k,n,m\ge 1} \frac{\chi_d(k m)\lambda_f(nm)}{nmk} \int_{-\infty}^\infty W_1(t) W_2(t) V_1^+(k n,t)V_2(k m,t)V_3(nm,t) \frac{d^*t}{2\pi^2},
\end{align}
an off-diagonal part with Kloosterman sums,
\begin{multline}
\label{eq:KuzKloos}
 \frac{4}{T^2} \sum_{\pm} \sum_{k,n,m,r\ge 1} \frac{\chi_d(k m)\lambda_f(r)}{k \sqrt{nmr}} \sum_{c \ge 1} \frac{S(\pm nm,r, c)}{c}\\
 \times \int_{-\infty}^\infty \mathcal{J}^\pm(\tfrac{\sqrt{nmr}}{c}, t) W_1(t) W_2(t) V_1^\pm(k n,t)V_2(k m,t)V_3(r,t) \frac{d^* t}{2 \pi},
\end{multline}
and an Eisenstein series part,
\begin{align}
\label{eq:KuzEisen}
 \frac{4}{T^2} \sum_{\pm} \sum_{k,n,m,r\ge 1} \frac{\chi_d(k m)\lambda_f(r)}{k \sqrt{nmr}} \int_{-\infty}^\infty \frac{\lambda(nm, t) \lambda(r, -t)}{|\zeta(1+2it)|^2} W_1(t) W_2(t) V_1^\pm(k n,t)V_2(k m,t)V_3(r,t) \frac{dt}{2\pi}.
\end{align}
The diagonal part gives the main contribution, and we will bound the other two parts by a negative power of $T$.

\section{The Eisenstein series}
Reversing the step where we combined the Hecke eigenvalues at $m$ and $n$, the Eisenstein series part \eqref{eq:KuzEisen} is
\begin{align}
	& \frac{4}{T^2} \sum_\pm \sum_{n,m,r\ge 1} \int_{-\infty}^\infty \frac{\lambda(n, t) \lambda(m, t) \lambda(r, -t) \chi_d(m) \lambda_f(r)}{|\zeta(1+2it)|^2 \sqrt{nmr}} W_1(t) W_2(t) V_1^\pm(n,t)V_2(m,t)V_3(r,t) \frac{dt}{2\pi}.
\end{align}
By Lemma \ref{afelem2}, this equals
\begin{align}
\label{eisin}	 \frac{4}{T^2} \int_{-\infty}^\infty W_1(t) W_2(t) \frac{|\zeta(\frac{1}{2}+it)|^2 |L(\frac{1}{2}+it,\chi_d)|^2 |L(\frac{1}{2}+it,f)|^2}{|\zeta(1+2it)|^2} \frac{dt}{2\pi}.
\end{align}
Now applying the bound (see \cite[chapter 5]{iwakow}),
\begin{align}
\zeta(1+2it) \gg \log(1+|t|)^{-1},
\end{align} 
and the subconvex bounds (see \cite[chapter 8]{iwakow}),
\begin{align}
\zeta(\thalf + it) \ll (1+|t|)^{\frac{1}{6}+\epsilon}, \quad L(\thalf + it, \chi_d) \ll (1+|t|)^{\frac{1}{6}+\epsilon},
\end{align}
it follows that (\ref{eisin}) is bounded by
 \begin{align}
 \label{eisin2} \frac{1}{T^{\frac{4}{3}-\epsilon}} \int_{0}^{T^{1+\epsilon}}  |L(\tfrac{1}{2}+it,f)|^2 dt.
\end{align}
Note that $L(\frac{1}{2}+it,f)$ has analytic conductor of size $1+|T^2-t^2|$, so that in the integral above, we may replace it by an approximate functional equation of length about $T$. Now by \cite[theorem 9.1]{iwakow}, we see that (\ref{eisin2}) is less than $T^{-\frac{1}{3}+\epsilon}$. 

\section{The diagonal} \label{diagsection}
Using (\ref{hmult}) and (\ref{v1approx}-\ref{v2approx}), the diagonal \eqref{eq:KuzDiag} equals
\begin{multline}
 \frac{4}{T^2}  \int_{-\infty}^\infty W_1(t) W_2(t) \frac{1}{(2\pi i )^3} \int_{(\epsilon)} \int_{(\epsilon)} \int_{(\epsilon)} e^{s_1^2+s_2^2+s_3^2} T^{- \beta s_1} d^{s_2} \left(\frac{|t|}{2\pi }\right)^{s_1+s_2} \left( \frac{|T^2-t^2|}{4 \pi^2 } \right)^{s_3} \zeta(1+2s_3)   \\
	\times \frac{L(1+s_1+s_2,\chi_d) L(1+s_1+s_3,f) L(1+s_2+s_3,f \times \chi_d)}{L(2+s_1+s_2+2s_3,\chi_d)} \frac{ds_1 \, ds_2 \, ds_3}{s_1 s_2 s_3} \frac{d^*t}{2\pi^2} +O(T^{-1+\alpha+\epsilon}).
\end{multline}
We shift each line of integration back to $\Re(s_i) = -\frac{1}{2}+\epsilon$, picking up poles at $s_i=0$. By the rapid decay of $e^{s_i^2}$, we may restrict each integral to $|\Im(s_i)|<T^{\epsilon}$.
The residue at $s_1=s_2=s_3=0$ gives the main term, and the three terms involving integrals over the shifted contours may be bounded by a negative power of $T$ using subconvexity results as follows.
Consider the result of the first shift to $\Re(s_1) = -\frac{1}{2}+\epsilon$. Trivially bounding the short $s_i$-integrals, it is sufficient to bound
\begin{align}
\label{shiftedint} & T^{\frac{\beta}{2}-2+\epsilon} \sup_{y_1,y_2,y_3} \left| L(\tfrac{1}{2}+\epsilon+iy_1+iy_2,f) \right| \int_{0}^{T^{1+\epsilon}} t^{-\half} \ d^*t,
\end{align}
where we have set $y_i = \Im(s_i)$, and the supremum is taken over all $y_i \in (-T^\epsilon,T^\epsilon)$.
Since the $t$-integral is bounded by $T^{\frac{3}{2}+\epsilon}$, we have by the subconvex estimate \eqref{eq:fSubconv}, which holds a fortiori to the right of $\frac{1}{2}$, that (\ref{shiftedint}) is less than a negative power of $T$ if $\beta < 2\Delta$.
The other two error terms are similar. The first two residues in $s_1$ and $s_2$ are from the simple poles of $(s_1 s_2)^{-1}$, but the final residue at $s_3=0$ requires some additional work due to the double pole of $s_3^{-1} \zeta(1+2s_3) \sim \frac{1}{2s_3^2}$. The result of the final shift is
\begin{align}
 	\frac{4}{T^2} \int_0^\infty W_1(t) W_2(t) \mathop{\mathrm{res}}_{s_3=0} \Biggl(\frac{e^{s_3^2}}{s_3^2} \frac{L(1,\chi_d) L(1+s_3,f) L(1+s_3,f \times \chi_d)}{L(2+2s_3,\chi_d)} \left( \frac{|T^2-t^2|}{4 \pi^2 } \right)^{s_3} \Biggr) \frac{d^*t}{2\pi^2}.
\end{align}
The residue evaluates to
\begin{align}
	 L'(0) - 2L(0)\log 2\pi + L(0)\log |T^2-t^2|,
\end{align}
where we define
\begin{align}
L(s) = \frac{L(1,\chi_d) L(1+s,f) L(1+s,f \times \chi_d)}{L(2+2s,\chi_d)}.
\end{align}

Thus we arrive at the main term,
\begin{align}
\label{diag01} & \frac{4}{T^2} \int_0^\infty W_1(t) W_2(t) \left(  L'(0) - 2L(0)\log 2\pi + L(0)\log |T^2-t^2| \right) \frac{d^*t}{2\pi^2}.
\end{align}
Recall that on the intervals $(0,T^{1-\alpha})\cup(T-T^{1-\alpha},T+T^{1-\alpha}),$
the function $W_2(t)$ is negligible.
On the remaining ranges, the hyperbolic tangent in $d^*t$ may be replaced by 1 up to admissible error. Making the substitution $t\mapsto T\sqrt{t}$, we get that (\ref{diag01}) equals
\begin{align}
\label{eq:MainBeforeTransitions}
 & \frac{1}{\pi^2} \int_{\mathcal{I}_1} e^{-t} W_2(T \sqrt{t}) \left( 2L(0) \log T + L'(0) - 2L(0)\log 2\pi + L(0)\log |1-t|  \right) dt  +O(T^{-\alpha+\epsilon}),
\end{align}
where
\begin{align}
\mathcal{I}_1 = (T^{-\alpha},1-T^{-\alpha}) \cup (1+T^{-\alpha},\infty). 
\end{align}
On the transitional intervals $(T^{-\alpha},T^{-\frac{\alpha}{4}})\cup(1-T^{-\frac{\alpha}{4}},1-T^{-\alpha})\cup(1+T^{-\alpha},1+T^{-\frac{\alpha}{4}}),$
we apply the bounds $0<W_2(T \sqrt{t})<1$, and outside this range, $W_2$ is very close to 1. So \eqref{eq:MainBeforeTransitions} equals
\begin{align}
\label{diag02}& \frac{1}{\pi^2} \int_{\mathcal{I}_2} e^{-t} \left(2L(0) \log T + L'(0) - 2L(0)\log 2\pi + L(0)\log |1-t|  \right) dt+O(T^{-\frac{\alpha}{4}+\epsilon})\\
 & = \frac{2L(0) \log T + L'(0) - 2L(0)\log 2\pi}{\pi^2}   + \frac{L(0)}{\pi^2} \int_{\mathcal{I}_2} e^{-t} \log |1-t|  dt+O(T^{-\frac{\alpha}{4}+\epsilon}),
\end{align}
where
\begin{align}
\mathcal{I}_2 = (T^{-\frac{\alpha}{4}},1-T^{-\frac{\alpha}{4}}) \cup (1+T^{-\frac{\alpha}{4}},\infty).
\end{align}
This last integral may be evaluated using the exponential integral function \cite[Ch 5, see p228, footnote 3]{AbramStegun},
\begin{align}
 \mathrm{Ei}(x)=-\text{PV}\int_{-x}^\infty e^{-t} \frac{dt}{t} = \int_0^x \frac{e^t-1}{t} \, dt+\log |x|+\gamma, 
\end{align}
for $x$ real, where $\gamma$ is Euler's constant. We have that
\begin{align}
 \int_{\mathcal{I}_2} e^{-t} \log |1-t|  dt = \frac{1}{e}  \left(-\mathrm{Ei}(1)+\mathrm{Ei}(T^{-\frac{\alpha}{4}})-\mathrm{Ei}(-T^{-\frac{\alpha}{4}})\right)  +O(T^{-\frac{\alpha}{4}+\epsilon}).
\end{align}
Now from the second integral representation above, the exponential integral function satisfies the asymptotic
\begin{align}
 \mathrm{Ei}(x)=\log|x|+\gamma+O(x)
\end{align}
for $|x| < 1$, and this gives the main term in Theorem \ref{main} with 
\begin{align}
C= \frac{ L'(0) - 2L(0)\log 2\pi}{\pi^2}   - \frac{L(0)\mathrm{Ei}(1)}{\pi^2 e}.
\end{align}

\section{The off-diagonal: shorter sum}

We consider the shorter sum in the off-diagonal (\ref{eq:KuzKloos}):
\begin{multline}
 \frac{4}{T^2}  \sum_{k,n,m,r\ge 1} \frac{\chi_d(k m)\lambda_f(r)}{k \sqrt{nmr}} \sum_{c \ge 1} \frac{S(\pm nm,r, c)}{c}\\
 \times \int_{-\infty}^\infty \mathcal{J}^+\Big(\tfrac{\sqrt{nmr}}{c}, t \Big) W_1(t) W_2(t) V_1^\pm(k n,t)V_2(k m,t)V_3(r,t) \frac{d^* t}{2 \pi}.
\end{multline}
It suffices to restrict the $t$-integral above to (\ref{wrange}). It also suffices to treat the leading terms of (\ref{v1approx0}-\ref{v3approx0}), as the lower order terms are similar, and the part of the sum with $k$=1, as the terms with $k>1$ are similar. Thus we must bound by a negative power of $T$ the sum
\begin{align}
\label{shortsum} \frac{1}{T^2} \sum_{n,m,r,c\ge 1} \frac{\lam(r) \chi_d(m)}{\sqrt{nmr}}  \frac{S( nm, r, c)}{c} \intt \frac{J_{2it} \big(\frac{ 4\pi \sqrt{nmr}}{c}  \big)}{\cosh(\pi t)}  Z\left( \frac{t}{T} \right)  V^+ \left(\frac{n}{T} ,\frac{m}{T} , \frac{r}{T^2} ; \frac{t}{T} \right) tdt.
\end{align}
We apply Lemma \ref{jbeslemma}. In this application The main term of (\ref{jbeslem}) is less than $T^{-100}$ unless
\begin{align}
\label{cbound} c<  \frac{\sqrt{nmr}}{T^{2}}T^{{3\alpha}} \ll \frac{\sqrt{T^{1-\beta} T^{3+\epsilon}}}{T^{2}}T^{{3\alpha}}.
\end{align}
We now fix $\beta = 7 \alpha$, so that (\ref{cbound}) is impossible for a positive integer $c$ when $T$ is large enough. The error term $O(T^{-100})$ of Lemma \ref{jbeslemma} is dominated by the error term $O(\frac{\sqrt{nmr}}{cT^{3-12\alpha}})$ once any reasonable bound on $c$ is imposed, such as $c\le T^{10}$. This can be achieved by shifting the line of integration in (\ref{shortsum}) to $\Im(t)=-\frac{1}{2} +\epsilon$ and bounding absolutely to see that the contribution of larger $c$ is negligible. Thus only the error term $O(\frac{\sqrt{nmr}}{cT^{3-12\alpha}})$ of Lemma \ref{jbeslemma} contributes to (\ref{shortsum}), and this contribution is bounded by
\begin{align}
\frac{1}{T^2} \sum_{\substack{c\ge 1\\ n,m \le T^{1+\epsilon}\\r\le T^{2+\epsilon}}} \frac{1}{\sqrt{nmr}}  \frac{|S(nm, r, c)|}{c} \frac{\sqrt{nmr}}{cT^{3-12\alpha}} \ll T^{-1+12\alpha+\epsilon},
\end{align}
on using Weil's bound for the Kloosterman sum. 
\section{The off-diagonal: longer sum}

We now consider the longer sum in the off-diagonal, which we will bound by $T^{-\delta}$ for some absolute constant $\delta>0$, as long as $\alpha$ is small enough. Since the actual size of $\alpha$ does not affect the final bound, it will be very convenient to rename $\alpha$ to $\epsilon$ for this section, in order to employ the $\epsilon$-convention. We must bound by a negative power of $T$ the sum  
\begin{multline}
 \frac{1}{T^2} \sum_{n,m,r,c\ge 1} \frac{\lam(r) \chi_d(m)}{\sqrt{nmr}}  \frac{S(- nm, r, c)}{c} \\ \times \intt \sinh(\pi t) K_{2it} \Big( \frac{4\pi \sqrt{nmr}}{c} \Big) Z\left( \frac{t}{T} \right)  V^- \left(\frac{n}{T} ,\frac{m}{T} , \frac{r}{T^2} ; \frac{t}{T} \right) tdt.
 \end{multline}
We apply Lemma \ref{kbeslemma}. As before, the error term $O(T^{-100})$ can be ignored. The contribution of the error term $O(\frac{\sqrt{nmr}}{cT^{4-\epsilon}})$ from (\ref{kbeslem}) is bounded by
\begin{align}
\frac{1}{T^2} \sum_{\substack{c\ge 1\\ n,m \le T^{1+\epsilon}\\r\le T^{2+\epsilon}}} \frac{1}{\sqrt{nmr}}  \frac{|S(- nm, r, c)|}{c} \frac{\sqrt{nmr}}{cT^{4-\epsilon}} \ll T^{-2+\epsilon},
\end{align}
on using Weil's bound for the Kloosterman sum. Thus it suffices to consider only the main terms of (\ref{kbeslem}). The second main term is a non-oscillatory bump function like the first, but of lower order. Thus it suffices to treat only the leading term of (\ref{kbeslem}). We must bound by a negative power of $T$ the sum
\begin{multline}
 \frac{1}{T^2} \sum_{n,m,r,c\ge 1} \frac{\lam(r) \chi_d(m)S(-nm, r, c)}{c^2}  Z\left( \frac{2\pi \sqrt{nmr}}{T c} \right)  V^- \left(\frac{n}{T} ,\frac{m}{T} , \frac{r}{T^2} ;  \frac{2 \pi \sqrt{nmr}}{T c} \right).
\end{multline}
Applying a smooth partition of unity, we consider the sum above in dyadic intervals. For $U$ a smooth bump function supported on $(1,2)\times(1,2)\times(1,2)$, it suffices to bound by a negative power of $T$ the sum
\begin{multline}
\label{dyadic} \frac{1}{T^2} \sum_{n,m,r,c\ge 1} \frac{\lam(r) \chi_d(m)S(-nm, r, c)}{c^2} 
Z\left( \frac{2\pi \sqrt{nmr}}{T c} \right)  V^- \left(\frac{n}{T} ,\frac{m}{T} , \frac{r}{T^2} ;  \frac{2\pi \sqrt{nmr}}{T c} \right) U\left( \frac{n}{N} ,\frac{m}{M},\frac{r}{R}\right)
\end{multline}
for
\begin{align}
\label{ranges} N<T^{1+\epsilon}, \ M<T^{1+\epsilon} \text{ and } R< T^{2+\epsilon}.
\end{align}
The function $Z$ restricts the sum to 
\begin{align}
\label{ranges2} T^{-1}\sqrt{NMR} \ll c \ll T^{-1+\epsilon} \sqrt{NMR}.
\end{align}

\subsubsection{{\bf Poisson summation}}

\

{\bf Case I.} Suppose that $d|c$. Then we replace $c$ by $cd$ in (\ref{dyadic}) and apply Poisson summation in $n$ and $m$ (after splitting into residue classes modulo $cd$) to get that the part of (\ref{dyadic}) with $d|c$ is bounded by
\begin{align}
\label{afterpoiss} \frac{NM}{T^2}  \sum_{-\infty<\ell_1,\ell_2<\infty} \ \sum_{ r,c\ge 1} \frac{\lam(r)}{c^4}  \sum_{a_1,a_2 \bmod cd } \chi_d(a_2) S(-a_1a_2, r, cd) e\left(\frac{a_1\ell_1+a_2\ell_2}{cd}\right)  \phi \left( \frac{r}{R}, \frac{cT}{\sqrt{NMR}} \right),
\end{align}
where
\begin{multline}
\label{fdef}  \phi_{\ell_1,\ell_2}\left( y_1, y_2 \right) = \phi\left( y_1, y_2 \right)=  \intt \intt Z\left( \frac{2\pi \sqrt{x_1 x_2 y_1}}{y_2d} \right) V^-\left( \frac{x_1 N}{T},\frac{x_2M}{T} , \frac{y_1R}{T^2}; \frac{2\pi \sqrt{x_1 x_2 y_1}}{y_2d} \right) \\ \times U\left(x_1 , x_2  ,y_1 \right)  e\left( \frac{-N \ell_1 x_1-M\ell_2 x_2}{y_2 dT^{-1} \sqrt{NMR}} \right) \ dx_1dx_2.
\end{multline}
Writing 
\begin{align}
S(-a_1a_2, r, cd) = \sums_{a_3 \bmod cd} e\left(\frac{-a_1a_2\overline{a_3}+ra_3}{cd}\right),
\end{align}
we have that
\begin{align}
\label{expev} \sum_{a_1,a_2 \bmod cd } \chi_d(a_2) S(-a_1a_2, r, cd) e\left(\frac{a_1\ell_1+a_2\ell_2}{cd}\right) = cd \sums_{a_3 \bmod cd} \chi_d(\ell_1 a_3)e\left(\frac{ a_3(r+ \ell_1 \ell_2 )}{cd}\right).
\end{align}
So (\ref{afterpoiss}) is bounded by
\begin{align}
\label{afterpoi} \frac{NM}{T^2} \sum_{-\infty<\ell_1,\ell_2<\infty} \ \sum_{ r,c\ge 1} \frac{\lam(r)}{c^3} \sums_{a_3 \bmod cd} \chi_d(\ell_1 a_3)e\left(\frac{ a_3(r+ \ell_1 \ell_2 )}{cd}\right)  \phi\left( \frac{r}{R}, \frac{cT}{\sqrt{NMR}} \right).
\end{align}

By repeatedly integrating by parts the integral in (\ref{fdef}), we see that $  \phi \left( \frac{r}{R}, \frac{cT}{\sqrt{NMR}} \right)  \ll T^{-100}$ unless
\begin{align}
\label{lsize1} |\ell_1| < \frac{c}{N}T^{\epsilon} \text{ and } |\ell_2| < \frac{c}{M}T^{\epsilon}.
\end{align}
Thus by (\ref{ranges}) and (\ref{ranges2}) we may assume that
\begin{align}
\label{lsize2} |\ell_1 \ell_2| \ll \frac{c^2}{NM}T^{\epsilon} \ll T^{\epsilon}.
\end{align}
This implies that if $\ell_1$ and $\ell_2$ are non-zero then they are both less than $T^{\epsilon}$, so that our notation is suggestive in suppressing the dependence of $\phi$ on $\ell_1$ and $\ell_2$. If $\ell_1=0$ then $\chi_d(\ell_1)=0$ and (\ref{afterpoi}) vanishes. In section \ref{zerol} we show that the contribution to (\ref{afterpoi}) of the terms with $\ell_2=0$ is small. Hence we assume that $\ell_1\ell_2\neq 0$, so that it suffices to bound by a negative power of $T$ the sum
\begin{align}
\label{b4vor} \frac{NM}{T^2} \sum_{ r,c\ge 1} \frac{\lam(r)}{c^3} \sums_{a_3 \bmod cd} \chi_d( a_3)e\left(\frac{ a_3 r}{cd}\right)  e\left(\frac{ a_3\ell_1 \ell_2 }{cd}\right) \phi \left( \frac{r}{R}, \frac{cT}{\sqrt{NMR}} \right)
 \end{align}
for any $0< |\ell_1|,|\ell_2|< T^{\epsilon} $. We must have by (\ref{lsize1}) that $N<cT^{\epsilon}$ and $M<cT^{\epsilon}$. This implies by (\ref{ranges}) and (\ref{ranges2}) that we must have
\begin{align}
\label{newranges} T^{-\epsilon} < \frac{N}{M} < T^{\epsilon}, \ \ \ \ T^{-\epsilon} < \frac{R}{T^2} < T^{\epsilon}, \ \ \ \ T^{-\epsilon} < \frac{c}{N} < T^{\epsilon}.
\end{align}
Note that in these ranges we have
\begin{align}
\label{goodderivs} \frac{\partial^{n_1+n_2}}{\partial y_1^{n_1}\partial y_2^{n_2}} \phi(y_1,y_2) \ll T^{\epsilon(n_1+n_2)}.
\end{align}
Here we used (\ref{lsize2}) and (\ref{newranges}) to see that in the exponential factor $e\left( \frac{-N \ell_1 x_1-M\ell_2 x_2}{y_2d T^{-1} \sqrt{NMR}} \right)$ of (\ref{fdef}), we have $ \frac{N \ell_1 x_1}{d T^{-1} \sqrt{NMR}} < T^{\epsilon}$ and $\frac{M \ell_2 x_1}{d T^{-1} \sqrt{NMR}} < T^{\epsilon}$.

{\bf Case II.} Suppose that $(c,d)=1$. Then by Poisson summation in $n$ (after splitting into residue classes modulo $c$) and $m$ (after splitting into residue classes modulo $cd$), we get that the part of (\ref{dyadic}) with $(c,d)=1$ is bounded by

\begin{align}
\label{varafterpoiss} \frac{NM}{T^2}  \sum_{-\infty<\ell_1,\ell_2<\infty} \ \sum_{ \substack{r,c\ge 1\\ (c,d)=1}} \frac{\lam(r)}{c^4}  \sum_{\substack{a_1 \bmod c \\ a_2 \bmod cd } } \chi_d(a_2) S(-a_1a_2, r, c) e\left(\frac{a_1d\ell_1+a_2\ell_2}{cd}\right)  \varphi \left( \frac{r}{R}, \frac{cT}{\sqrt{NMR}} \right),
\end{align}
where
\begin{multline}
\label{varfdef}  \varphi_{\ell_1,\ell_2}\left( y_1, y_2 \right)=\varphi\left( y_1, y_2 \right)=  \intt \intt Z\left( \frac{2\pi \sqrt{x_1 x_2 y_1}}{y_2} \right) V^-\left( \frac{x_1 N}{T},\frac{x_2M}{T} , \frac{y_1R}{T^2}; \frac{2\pi \sqrt{x_1 x_2 y_1}}{y_2} \right) \\ \times  U\left(x_1 , x_2  ,y_1 \right) e\left( \frac{-N d\ell_1 x_1-M\ell_2 x_2}{y_2 dT^{-1} \sqrt{NMR}} \right) \ dx_1dx_2.
\end{multline}
By the Chinese Remainder Theorem, we have that
\begin{align}
\sum_{a_2 \bmod cd } \chi_d(a_2) S(-a_1a_2, r, c) e\left(\frac{a_2\ell_2}{cd}\right) 
= \sum_{a_2 \bmod c } S(-a_1a_2 d, r, c) e\left(\frac{a_2\ell_2}{c}\right) \sum_{b_2 \bmod d } \chi_d(b_2 c) e\left(\frac{b_2\ell_2}{d}\right).
\end{align}
The innermost sum equals $ \chi_d(c\ell_2) d^\half$. Thus, as in (\ref{expev}), we get that (\ref{varafterpoiss}) is bounded by
\begin{align}
\label{b4vorvar} \frac{NM}{T^2} \sum_{-\infty<\ell_1,\ell_2<\infty} \ \sum_{ r,c\ge 1} \frac{\lam(r)}{c^3} \sums_{a_3 \bmod c} \chi_d(c \ell_2)e\left(\frac{ a_3(r+ \overline{d} \ell_1 \ell_2 )}{c}\right)  \varphi \left( \frac{r}{R}, \frac{cT}{\sqrt{NMR}} \right).
\end{align}
Statements analogous to (\ref{lsize1}-\ref{goodderivs}) hold for the sum above.

 \subsubsection{{\bf Voronoi summation}} \label{vorsect}
 
 By Lemma \ref{voronoilem} we have that (\ref{b4vor}) equals
 \begin{align}
\label{vorstate} \frac{NM}{T^2} \sum_{ q,c\ge 1} \frac{\lam(q)}{ q c^2} \sum_\pm S_{\chi_d} (\pm q , \ell_1\ell_2, cd)   \int_{(\sigma)} \int_0^\infty \phi\left( x, \frac{cT}{\sqrt{NMR}} \right)x^{-s-1} dx  \left(\frac{\pi^2 Rq}{c^2d^2}\right)^{-s} G_\pm(s) ds.
 \end{align}
Writing $s=\sigma+i\gamma$, observe that by (\ref{goodderivs}), we may restrict the $s$-integral in (\ref{vorstate}) to $|\gamma|<T^{\epsilon}$, else the integral is less than $T^{-100}$ by integration by parts. In this range, by Stirling's approximation (as in \ref{v1approx0}-\ref{v3approx0})) we have that
\begin{align}
\label{gplus} & G_+(s) = \frac{-i}{2\pi^2} \left( \frac{T}{2} \right)^{2s+1}  \Big( 1 + \sum_{n=1}^{1000} \frac{C_n(\sigma,\gamma)}{T^{n}} +O(T^{-100})\Big),\\
 &G_-(s) = \frac{-i}{2\pi^2} \left( \frac{T}{2} \right)^{2s+1}  \Big(  \sum_{n=1}^{1000} \frac{C_n(\sigma,\gamma)}{T^{n}} +O(T^{-100})\Big),
\end{align}
for some $C_n(\sigma,\gamma)$ (not necessarily the same in each expression above) polynomial in $\sigma$ and $\gamma$.
With this observation, by moving the integral in (\ref{vorstate}) far to the right or left, without crossing any poles as we restrict to $|\gamma|<T^\epsilon$, we find that we can restrict the $q$-sum to $ \frac{c^2T^{2-\epsilon}}{R} < q < \frac{c^2T^{2+\epsilon}}{R} $, or equivalently by (\ref{newranges}), to 
\begin{align}
\label{qrange} c^2T^{-\epsilon}< q < c^2T^{\epsilon}.
\end{align}
 Thus to treat (\ref{b4vor}) after Voronoi summation, it suffices to consider only the leading term of (\ref{gplus}) as the rest are similar, and to bound by a negative power of $T$ the sum
\begin{align}
 \frac{NM}{T} \sum_{ q,c\ge 1} \frac{\lam(q)}{ q  c}  S_{\chi_d} ( q , \ell_1\ell_2, cd)  \Phi\left( \frac{4\pi \sqrt{q|\ell_1 \ell_2|}}{cd}  \right),
\end{align}
where
 \begin{align}
 \label{ffdef} \Phi(y) = X(y)\int_{(1)}  \int_0^\infty \phi\left( x, \frac{4\pi T\sqrt{q|\ell_1\ell_2|}}{yd\sqrt{NMR}} \right) x^{-s-1}  dx \left( \frac{  y  \sqrt{R}}{2 T\sqrt{|\ell_1\ell_2|}} \right)^{-2s} ds
 \end{align}
 and $X$ is any smooth compactly supported function on $(T^{-\epsilon}, T^{\epsilon})$ whose derivatives are bounded by powers of $T^{\epsilon}$. By (\ref{newranges}) and (\ref{qrange}), this is equivalent to bounding by a negative power of $T$ the sum 
\begin{align}
\label{tobound} \frac{1}{T} \sum_{ q} \frac{\lam(q)}{ q^\half  }  \sum_{c\ge 1} \frac{S_{\chi_d} ( q , \ell_1\ell_2, cd)}{c}  \Phi\left( \frac{4\pi \sqrt{q|\ell_1 \ell_2|}}{cd}  \right).
\end{align}
 By (\ref{newranges}) and (\ref{qrange}), we have that
 \begin{align}
 \label{nicederivs} \Phi^{(n)}(y) \ll (T^{\epsilon})^n.
 \end{align}
 
When Voronoi summation is applied to (\ref{b4vorvar}), we obtain the exponential sum $ \chi_d(c\ell_2)S(q,\overline{d}\ell_1\ell_2,c)$. We may assume that $(\ell_2,d)=1$ or else this vanishes. Writing $\ell_1=d^k \ell_3$, where $k\ge 0$ and $(\ell_3,d)=1$, we have, by a multiplicative property of Kloosterman sums \cite[Section 2]{blomil} for $(c,d)=1$, that
\begin{align}
\chi_d(c\ell_2)S(q,\overline{d}\ell_1\ell_2,c) &= \chi_d(c\ell_2)S(d^k q,\overline{d}\ell_3\ell_2,c)\\
\nonumber &=  \chi_d(c\ell_2) S_{\chi_d}(0,\overline{c} \ell_3\ell_2,d)^{-1} S_{\chi_d}(d^{k+1}q,\ell_3\ell_2,cd)\\
\nonumber &= d^{-\half} \chi_d(\ell_3 )  S_{\chi_d}(d^{k+1}q,\ell_3\ell_2,cd).
\end{align}
Further note that the condition $(c,d)=1$ may be dropped because when $d|c$,
\begin{align}
\label{klozero} S_{\chi_d}(d^{k+1}q,\ell_3\ell_2,cd)=0
\end{align}
by the argument in the next paragraph. Thus in this case Voronoi summation leads to an expression similar to (\ref{tobound}). We therefore show the details of the rest of the proof for only (\ref{tobound}).

Before moving on, we prove (\ref{klozero}). Write $c=c'd^{1+j}$, where $j\ge 0$ and $(c',d)=1$. Then the left hand side of (\ref{klozero}) is a multiple of 
\begin{align}
\label{klo0} S_{\chi_d}(d^{k+1}q\overline{c'},\ell_3\ell_2\overline{c'}, d^{2+j})= \sums_{a  \bmod d^{2+j} } \chi_d(a) e\Big( \frac{d^{k+1} q\overline{ac'}+\ell_3\ell_2 \overline{c'} a}{d^{2+j}} \Big).
\end{align}
We may write $a=v+ud^{1+j}$, where $u$ ranges over all residue classes modulo $d$ and $v$ ranges over the primitive residue classes modulo $d^{1+j}$. Note that $d(\overline{v+ud^{1+j}}) \equiv d \overline{v}$ modulo $d^{2+j}$, so that (\ref{klo0}) equals
\begin{align}
\sum_{u \bmod d} \  \sums_{v \bmod d^{j+1}} \chi_d(v) e\Big( \frac{d^{k+1} q \overline{c' v} + \ell_3\ell_2 \overline{c'} (v+ud^{1+j})}{d^{2+j}}\Big).
\end{align}
The $u$-sum vanishes as $(\overline{c'}\ell_3 \ell_2,d)=1$.

\subsubsection{{\bf Kuznetsov's formula and subconvexity}}

Applying Lemma \ref{kuzbackwards} and the remarks following it to the $c$-sum of (\ref{tobound}), we find that it suffices to bound by a negative power of $T$ the sums
\begin{align}
\label{fnlbnd1} \frac{1}{T} \sum_{q\ge 1} \frac{\lam(q) \lambda_g(q)}{ q^\half } \phi\left( x, \frac{4\pi T\sqrt{q|\ell_1\ell_2|}}{yd\sqrt{NMR}} \right) 
\end{align}
and
\begin{align}
\label{fnlbnd2} \frac{1}{T} \sum_{q\ge 1} \frac{\lam(q) \lambda_v(q,-t)}{ q^\half } \phi\left( x, \frac{4\pi T\sqrt{q|\ell_1\ell_2|}}{yd\sqrt{NMR}} \right),
\end{align}
where $g\in \mathcal{B}_k(d,\chi_d)$ with $k<T^\epsilon$ or $g\in \mathcal{B}(d,\chi_d)$ with $|t_g|<T^\epsilon$, $v|d$ and $|t|<T^{\epsilon}$, $x\in(1,2)$ and $y\in (T^{-\epsilon},T^{\epsilon})$.  By (\ref{newranges}) and (\ref{qrange}), the $q$-sums have length about $N$. If $N<T^{2-\delta}$ for some $\delta>0$ then (\ref{fnlbnd1}-\ref{fnlbnd2}) may be bounded trivially, while if $N>T^{2-\delta}$ for $\delta$ small enough, then the required bound follows by (\ref{subconvsum1}-\ref{subconvsum3})  and partial summation.

\subsubsection{{\bf The case $\boldsymbol{\ell_2=0}$}}\label{zerol}

When $\ell_1\neq 0$ and $\ell_2=0$,  we have that (\ref{afterpoi}) equals
\begin{align}
\label{l2} \frac{NM}{T^2}\sum_{-\infty<\ell_1<\infty} \ \sum_{ r,c\ge 1} \frac{\lam(r)}{c^3} \sums_{a \bmod cd} \chi_d(\ell_1 a)e\left(\frac{a r}{cd}\right)  \phi_{\ell_1,0} \left( \frac{r}{R}, \frac{cT}{\sqrt{NMR}} \right).
\end{align}
The first observation is that we may restrict (\ref{l2}) to $|\ell_1|<T^{\epsilon}$. To see this, let $T^{\epsilon} < L <\frac{c}{N}T^{\epsilon}$, by (\ref{lsize1}), and consider
\begin{align}
\label{l3} \frac{NM}{T^2}\sum_{-\infty<\ell_1<\infty} W\left(\frac{\ell_1}{L} \right) \sum_{ r,c\ge 1} \frac{\lam(r)}{c^3} \sums_{a \bmod cd} \chi_d(\ell_1 a)e\left(\frac{a r}{cd}\right)  \phi_{\ell_1,0} \left( \frac{r}{R}, \frac{cT}{\sqrt{NMR}} \right),
\end{align}
for any fixed smooth function $W$ compactly supported on $(1,2)$. By Poisson summation in $\ell_1$ (after splitting into residue classes modulo $d$), we get that (\ref{l2}) equals
\begin{multline}
\label{l4} \frac{NML}{T^2d}  \sum_{ r,c\ge 1} \frac{\lam(r)}{c^3} \sums_{a \bmod cd} \chi_d( a) e\left(\frac{a r}{cd}\right) \sum_{b \bmod d} \chi_d(b)  \\ \times  \sum_{-\infty<k<\infty}  e\left( \frac{bk}{d} \right)  \intt W(z)\phi_{zL,0}\left( \frac{r}{R}, \frac{cT}{\sqrt{NMR}} \right) e\left( \frac{-z L k}{d} \right)  dz.
\end{multline}
Since $L <\frac{c}{N}T^{\epsilon}$, we have that 
\begin{align}
\frac{\partial^{n}}{\partial z^{n}} W(z) \phi_{zL,0}(y_1,y_2) \ll (T^{\epsilon})^{n}.
\end{align}
Using this and that $L>T^{\epsilon}$, we find by repeatedly integrating by parts the $z$-integral in (\ref{l4}) that the contribution of $|k|\ge 1$ is less than $T^{-100}$, say. This leaves the contribution of $k=0$ to (\ref{l4}), which vanishes as $\sum_{b \bmod d} \chi_d(b)=0$.

It suffices now to bound by a negative power of $T$ the sum
\begin{align}
\label{l6} \frac{NM}{T^2}  \sum_{ r,c\ge 1} \frac{\lam(r)}{c^3} \sums_{a \bmod cd} \chi_d(a)e\left(\frac{a r}{cd}\right) \phi_{\ell_1,0} \left( \frac{r}{R}, \frac{cT}{\sqrt{NMR}} \right)
\end{align}
for any $|\ell_1| <T^{\epsilon}$. 
We split up this sum according to the value of $(c,d)=d^j$ for $j=0$ or $1$:
\begin{align}
\label{l7} \frac{NM}{T^2} \sum_{j=0,1} \sum_{\substack{r,c\ge 1\\(c,d^{1-j})=1}} \frac{\lam(r)}{c^3 d^{3j}} \sums_{\substack{a \bmod d^{j+1}\\ b \bmod c}} \chi_d(ac)  e\left(\frac{a r}{d^{j+1}}\right)  e\left(\frac{b r}{c}\right)  \phi_{\ell_1,0} \left( \frac{r}{R}, \frac{c d^j T}{\sqrt{NMR}} \right).
\end{align}
We show how to treat this in the case that $j=0$, the other case being similar. Using that $\displaystyle \sums_{a \bmod d} \chi_d(a) e\Big(\frac{ar}{d}\Big)= d^\half \chi_d(r) $, it suffices to bound by a negative power of $T$ the sum
\begin{align}
\label{l8} \frac{NM}{T^2}   \sum_{c\ge 1} \frac{1}{c^3} \left| \sum_{r\ge 1} \lam(r)\chi_d(r) \mathcal{R}_{c}(r)  \phi_{\ell_1,0} \left( \frac{r}{R}, \frac{c T}{\sqrt{NMR}} \right) \right|,
\end{align}
where
\begin{align}
 \mathcal{R}_{c}(r)={ \sums_{b \bmod c}} e\left(\frac{br}{c}\right)= \sum_{c'|(c,r)}\mu\left( \frac{c}{c'} \right) c'
\end{align}
is a Ramanujan sum. 
Using the above identity, it suffices to bound by a negative power of $T$ the sum
\begin{align}
\label{l9a} \frac{NM}{T^2}   \sum_{c\ge 1} \frac{1}{(c')^2 c^3} \left| \sum_{\substack{r\ge 1}} \lam(r c')\chi_d(r ) \phi_{\ell_1,0} \left( \frac{r c'}{R}, \frac{c c' T}{\sqrt{NMR}} \right) \right|
\end{align}
for any positive integer $c'$. By Hecke multiplicativity,
\begin{align}
\lam(rc') = \sum_{c''|(c',r)} \mu(c'') \lam\left(\frac{ r}{c''} \right) \lam\Big(\frac{c'}{c''} \Big),
\end{align}
it suffices to bound by a negative power of $T$ the sum
\begin{align}
\label{l9b} \frac{NM}{T^2}   \sum_{c\ge 1} \frac{1}{(c')^2(c''c)^3} \left| \sum_{\substack{r\ge 1}} \lam( r)\chi_d(r ) \phi_{\ell_1,0} \left( \frac{r c' c''}{R}, \frac{c c' c'' T}{\sqrt{NMR}} \right) \right|,
\end{align}
or equivalently,
\begin{align}
\label{l10} \frac{NM}{T^{2}}  \sum_{ \frac{\sqrt{NMR}}{c'c''T^{1+\epsilon}} \le  c\le \frac{\sqrt{NMR}}{c' c'' T^{1-\epsilon}}} \frac{1}{(c')^2(c''c)^3} \left| \sum_{\frac{R}{c' c''} < r < \frac{2R}{c' c''}}   \lam( r)\chi_d(r ) \phi_{\ell_1,0} \left( \frac{r c' c''}{R}, \frac{c c' c'' T}{\sqrt{NMR}} \right) \right|.\end{align}
for any positive integers $c', c''$. The trivial bound for this is easily verified to be $T^{\epsilon}$. To do better, note that
\begin{align}
\frac{R}{c' c''} \gg \frac{c^2 c' c'' T^{2-\epsilon}}{NM} \gg T^{1-\epsilon}.
\end{align}
The first inequality follows by (\ref{ranges2}) and the second holds because $c c' c'' > N T^{-\epsilon}$ by (\ref{lsize1}) as $|\ell_1|\ge 1$, and $M<T^{1+\epsilon}$ by (\ref{ranges}). Thus by (\ref{subconvsum2}) and partial summation, the required bound for (\ref{l10}) follows. 

\subsection*{ Acknowledgments.} We are grateful to Valentin Blomer and Matthew Young for very helpful discussions regarding this project. The first author was supported by a grant from the European Research Council (grant agreement number 258713) when he was based at the University of G\"{o}ttingen and thanks Texas A\&M University at Qatar, where part of this work was done, for its hospitality.

\bibliographystyle{amsplain}

\bibliography{tripleproduct}

\end{document}